# Using a New Nonlinear Gradient Method for Solving Large Scale Convex Optimization Problems with an Application on Arabic Medical Text

Jaafar Hammoud[a*], Ali Eisa[b], Natalia Dobrenko[a], Natalia Gusarova[a]

*[a]ITMO University, St. Petersburg, Russia*
*[*]hammoudgj@gmail.com*
*[b]Aleppo University, Aleppo, Syria*

**Abstract**

Gradient methods have applications in multiple fields, including signal processing, image processing, and dynamic systems. In this paper, we present a nonlinear gradient method for solving convex supra-quadratic functions by developing the search direction, that done by hybridizing between the two conjugate coefficients HRM [2] and NHS [1]. The numerical results proved the effectiveness of the presented method by applying it to solve standard problems and reaching the exact solution if the objective function is quadratic convex.

Also presented in this article, an application to the problem of named entities in the Arabic medical language, as it proved the stability of the proposed method and its efficiency in terms of execution time.

## 1. Introduction

Several nonlinear conjugate gradient methods have been presented for solving high-dimensional unconstrained optimization problems which are given as [3]:

$$\min f(x) \quad ; x \in R^n \tag{1}$$

To solve problem (1), we start with the following iterative relationship:

$$x_{(k+1)} = x_k + \alpha_k d_k, \quad k = 0,1,2, \tag{2}$$

where $\alpha_k > 0$ is a step size, that is calculated by strong Wolfe-Powell's conditions [4]

$$f(x_k + \alpha_k d_k) \leq f(x_k) + \delta \alpha_k g_k^T d_k,$$
$$|g(x_k + \alpha_k d_k)^T d_k| \leq \sigma |g_k^T d_k| \tag{3}$$

where $0 < \delta < \sigma < 1$

$d_k$ is the search direction that is computed as follow [3]

$$d_{k+1} = \begin{cases} -g_k & \text{,if } k = 0 \\ -g_{k+1} + \beta_k d_k & \text{,if } k \geq 1 \end{cases} \tag{4}$$

Where $g_k = g(x_k) = \nabla f(x_k)$ represent the gradient vector for $f(x)$ at the point $x_k$, $\beta_k \in \mathbb{R}$ is known as CG coefficient that characterizes different CG methods. Some classical methods such below:

$$\beta_k^{HS} = \frac{g_k^T(g_k - g_{k-1})}{(g_k - g_{k-1})^T d_{k-1}} \qquad \text{HS [4]}$$

$$\beta_k^{FR} = \frac{g_k^T g_k}{\|g_{k-1}\|^2} \qquad \text{FR [5]}$$

$$\beta_k^{PRP} = \frac{g_k^T(g_k - g_{k-1})}{\|g_{k-1}\|^2} \qquad \text{PRP [6, 7]}$$

$$\beta_k^{HRM} = \frac{g_k^T \left(g_k - \frac{\|g_k\|}{\|g_{k-1}\|} g_{k-1}\right)}{\tau \|g_{k-1}\|^2 + (1-\tau)\|d_{k-1}\|^2}; \qquad \tau = 0.4 \qquad \text{HRM [2]}$$

$$\beta_k^{NHS} = \frac{\|g_k\|^2 - \frac{\|g_k\|}{\|g_{k-1}\|} \max{0, g_k^T g_{k-1}}}{\max\{max\, 0, u g_k^T d_{k-1} + \|g_{k-1}\|^2, d_k^T y_{k-1}\}}; \quad u = 1.1 \qquad \text{NHS [1]}$$

The iterative solution stops when we reach a point $x_k$ where the condition $\|g_k\| \leq \varepsilon$ is fulfilled, where $\epsilon$ is a very small positive number. Among the most common methods that rely on the aforementioned strategy are Newton's methods [9], quasi-Newton methods [10, 11, 12], trust region methods [13, 14], and conjugated gradient methods [15, 16].

On the other hand, the optimization techniques and methods play one of the most important roles in training neural networks (NN), because it is used to reduce the losses by changing the attributes of NN such as weights and learning rate.

The effect of choosing one optimization algorithm over another has been studied previously by many researchers, and despite the continuous development of this aspect, the widespread platforms that are used in the field of machine learning and deep learning depend on a specific group of these algorithms such as ADAM [22], SGD with momentum [23], and RMSprop [24], but these platforms come with the possibility of creating our own optimizer.

In [25], the named entity problem was studied on Arabic medical text taken from three medical volumes issued by the Arab Medical Encyclopedia, the researchers used a BERT model [26] that was introduced by Google.

As an application of the presented method, we implement our optimizer on the same dataset and show the results of comparison with the previous one.

During the following sections, a hybrid method for solving Problem (1) will be presented, and then its convergence will be studied, the numerical results of the mentioned method will be presented, and in the end, an application in the field of Arabic medical text processing will prove the efficiency of the method.

## 2. The formula and its convergence

In the following, we show a nonlinear gradient method for solving convex functions with high dimensions by hybridizing two CG formulas [28], and the new formula is given as:

$$\beta_k^{AWHM} = (1 - \theta_k)\beta_k^{NHS} + \theta_k \beta_k^{HRM} \qquad (5)$$

From (5) we distinguish the following cases:

- Case 1: if $\theta_k = 0$ then $\beta_k^{AWHM} = \beta_k^{NHS}$.

- Case 2: if $0 < \theta_k < 1$ then we find the new value for $\theta_k$ by using the search direction that was introduced by [16] as below:

$$d_{k+1}^T y_k = -ts_k^T g_{k+1}; \quad t > 0 \tag{6}$$

Where is $y_k = g_{k+1} - g_k$ and $.s_k = x_{k+1} - x_k$.

The new search direction is given by relation:

$$d_{k+1} = -g_{k+1} + \beta_k^{AWHM} d_k \tag{7}$$

From (5) and (7) we find that:

$$d_{k+1} = -g_{k+1} + \left((1-\theta_k)\beta_k^{NHS} + \theta_k \beta_k^{HRM}\right) d_k \tag{8}$$

And from (8) and (6) we find that:

$$-g_{k+1}^T y_k + (1-\theta_k)\beta_k^{NHS} d_k^T y_k + \theta_k \beta_k^{HRM} d_k^T y_k = -ts_k^T g_{k+1}$$

$$\theta_k (\beta_k^{HRM} - \beta_k^{NHS}) d_k^T y_k = -ts_k^T g_{k+1} + g_{k+1}^T y_k - \beta_k^{NHS} d_k^T y_k$$

$$\theta_k^{new} = \frac{-ts_k^T g_{k+1} + g_{k+1}^T y_k - \beta_k^{NHS} d_k^T y_k}{(\beta_k^{HRM} - \beta_k^{NHS}) d_k^T y_k} \tag{9}$$

- Case 3: if $\theta_k = 1$ then $\beta_k^{AWHM} = \beta_k^{HRM}$.

*2.1 Algorithms steps*

Input: $x_0 \in \mathbb{R}^n$ a start point, $f$ a goal function, and $\varepsilon > 0$.

Step 0: calculate the gradient vector $g_0 = \nabla f(x_0)$, the initial search direction $d_0 = -g_0$, and the step size $\lambda_0 = 1/\|g_0\|$, then we put $k = 0$, if $\|g_k\| \leq \epsilon$ we stop, else we go to step 1.

Step 1: we calculate the new search direction that satisfied the strong Wolfe-Powell conditions.

$$f(x_k + \lambda_k d_k) - f_k \leq \delta \lambda_k g_k^T d_k$$
$$|g(x_k + \lambda_k d_k)^T d_k| \leq -\sigma g_k^T d_k$$

Where $\delta \in (0, 0.5)$, $\sigma \in (\delta, 1)$.

Step 2: set a new point $x_{k+1} = x_k + \lambda_k d_k$, if $\|g_k\| \leq \epsilon$ we stop else go to step 3.

Step 3: calculate , $s_k = x_{k+1} - x_k$, $.y_k = g_{k+1} - g_k$.

Step 4: find the value of each of $\beta_k^{NHS}, \beta_k^{HRM}$ from below relations:

$$\beta_k^{HRM} = \frac{g_k^T \left(g_k - \frac{\|g_k\|}{\|g_{k-1}\|} g_{k-1}\right)}{\tau \|g_{k-1}\|^2 + (1-\tau)\|d_{k-1}\|^2}; \qquad \tau = 0.4.$$

$$\beta_k^{NHS} = \frac{\|g_k\|^2 - \frac{\|g_k\|}{\|g_{k-1}\|} \max{0, g_k^T g_{k-1}}}{\max\{max\, 0, ug_k^T d_{k-1} + \|g_{k-1}\|^2, d_k^T y_{k-1}\}}; \qquad u = 1.1$$

Step 5: calculate $\theta_k^{new}$ by the relation below:

$$\theta_k^{new} = \frac{-ts_k^T g_{k+1} + g_{k+1}^T y_k - \beta_k^{NHS} d_k^T y_k}{(\beta_k^{HRM} - \beta_k^{NHS}) d_k^T y_k}$$

Step 6: if $0 < \theta_k^{new} < 1$ then calculate the $\beta_k^{AWHM}$ as below:

$$\beta_k = \beta_k^{AWHM} = (1 - \theta_k)\beta_k^{NHS} + \theta_k\beta_k^{HRM}$$

If $\theta_k^{new} = 0$ then $\beta_k = \beta_k^{NHS}$
If $\theta_k^{new} = 1$ then $\beta_k = \beta_k^{HRM}$

Step 7: set the new search direction with the relation:
$$d^{new} = -g_{k+1} + \beta_k d_k$$

Step 8: if $|g_{k+1}^T g_k| \geq 0.2\|g_{k+1}\|^2$ then $d_{k+1} = -g_{k+1}$ else $d_{k+1} = d^{new}$, after that find

$$\lambda_{k+1} = \lambda_k \times \frac{\|d_k\|}{\|d_{k+1}\|}$$

Step 9: set k = k + 1 and go to step 1.

*2.2 Convergence analysis*

The following assumptions are often used in previous studies of the conjugate gradient methods: [2, 18]

**Assumption A:**
$f(x)$ is bounded from below on the level set $\Omega = \{x \in \mathbb{R}^n \mid f(x) \leq f(x_0)\}$, where $x_0$ is the starting point.

**Assumption B:**
In some neighbourhood $N$ of $\Omega$, the objective function is continuously differentiable, and its gradient is Lipschitz continuous, that is, there exists a constant $L > 0$ such that.

$$\|g(x) - g(y)\| \leq l\|x - y\| \quad \forall \ x, y \in N$$

**Assumption C:**
$$\forall x \in \Omega \|g(x)\| \leq \Gamma; \quad \Gamma \geq 0$$

**Theorem 1:**
Suppose that the sequences $\{g_k\}$ and $\{d_k\}$ are generated by the presented method. Then the sequence $\{d_k\}$ possesses as the sufficient descent condition
$$g_k^T d_k \leq c\|g_k\|^2 \qquad \forall k \geq 0, \quad c > 0$$

**Proof:**
For $k = 0$ the relation (10) is fulfilled, because:
$$g_0^T d_0 = -\|g_0\|^2$$
For $k \geq 1$
$$d_{k+1} = -g_{k+1} + \beta_k^{AWHM} d_k$$

$$d_{k+1} = -g_{k+1} + \left((1 - \theta_k)\beta_k^{NHS} + \theta_k\beta_k^{HRM}\right) d_k$$

$$d_{k+1} = -(\theta_k g_{k+1} + (1 - \theta_k)g_{k+1}) + \left((1 - \theta_k)\beta_k^{NHS} + \theta_k\beta_k^{HRM}\right) d_k$$

$$d_{k+1} = \theta_k(-g_{k+1} + \beta_k^{HRM} d_k) + (1 - \theta_k)(-g_{k+1} + \beta_k^{NHC} d_k)$$

$$d_{k+1} = \theta_k d_{k+1}^{HRM} + (1 - \theta_k)d_{k+1}^{NHC}$$

We discuss according to the value of $\theta_k$ we find:

i. If $\theta_k = 0$ then $d_{k+1} = d_{k+1}^{NHC}$

$$g_{k+1}^T d_{k+1} = g_{k+1}^T d_{k+1}^{NHC} = g_{k+1}^T(-g_{k+1} + \beta_k^{NHC} d_k) \leq c_1 \|g_{k+1}\|^2$$

where $c_1 = \left(1 - \frac{1}{\mu}\right); \mu = 1.1$

ii. If $\theta_k = 1$ then $d_{k+1} = d_{k+1}^{HRM}$

$$g_{k+1}^T d_{k+1} = g_{k+1}^T d_{k+1}^{HRM} = g_{k+1}^T(-g_{k+1} + \beta_k^{HRM} d_k) \leq c_2 \|g_{k+1}\|^2$$

where $c_2 = \left(2 - \frac{1}{1-5\sigma}\right); \sigma = 0.001$

iii. If $0 < \theta_k < 1$ then we suppose that: $0 < m_1 \leq \theta_k \leq m_2 < 1$

$$g_{k+1}^T d_{k+1} = \theta_k g_{k+1}^T d_{k+1}^{HRM} + (1 - \theta_k) g_{k+1}^T d_{k+1}^{NHC}$$

$$g_{k+1}^T d_{k+1} \leq m_1 g_{k+1}^T d_{k+1}^{HRM} + (1 - m_2) g_{k+1}^T d_{k+1}^{NHC}$$

$$g_{k+1}^T d_{k+1} \leq m_1 c_2 \|g_{k+1}\|^2 + (1 - m_2) c_1 \|g_{k+1}\|^2$$

$$g_{k+1}^T d_{k+1} \leq c \|g_{k+1}\|^2, \qquad c = m_1 c_2 + (1 - m_2) c_1$$

**Theorem 2:**

If the assumptions are fulfilled, and since the search direction fulfills the condition of sufficient descent condition, then the presented method fulfills the property of global convergence. That is, if the following relationship is fulfilled:

$$\sum_{k \geq 1} \frac{1}{\|d_{k+1}\|^2} = \infty$$

Then

$$\lim_{k \to \infty} (\inf \|g_{k+1}\|) = 0$$

**Proof:**

If the gradient vector $g_k \neq 0$ then there is a constant $r > 0$ where $\|g_k\| > r \ \forall \ k \geq 0$, then from (5) and by using Lipschitz condition and Assumption C, we find:

$$|\beta_k^{AWHM}| \leq +|\beta_k^{NHC}| + |\beta_k^{HRM}|$$

According to [1]:

$$0 < \beta_k^{NHC} \leq \frac{g_{k+1}^T d_{k+1}}{g_k^T d_k}$$

And according to [2]:

$$0 < |\beta_k^{HRM}| \leq \frac{1}{2b}, \quad b = \frac{5\bar{\gamma}^2(\gamma + \bar{\gamma})}{2\gamma^3} > 1$$

$$0 < |\beta_k^{AWHM}| \leq \left|\frac{g_{k+1}^T d_{k+1}}{g_k^T d_k}\right| + \frac{1}{2b} \leq \frac{c_1 \|g_{k+1}\|^2}{c_3 \|g_k\|^2} + \frac{1}{2b} \leq \frac{c_1 \Gamma^2}{c_3 r^2} + \frac{1}{2b} = \omega$$

$$0 < |\beta_k^{AWHM}| \leq \omega$$

From the iterative relation:

$$x_{k+1} = x_k + \lambda_k d_k \Rightarrow x_{k+1} - x_k = \lambda_k d_k \Rightarrow s_k = \lambda_k d_k \Rightarrow d_k = \frac{s_k}{\lambda_k}$$

$$\lambda_k \geq \lambda^* > 0 \Rightarrow \frac{1}{\lambda_k} \leq \frac{1}{\lambda^*}$$

According to our search direction

$$d_{k+1} = -g_{k+1} + \beta_k^{AWHM} d_k$$

$$\|d_{k+1}\| = \|-g_{k+1} + \beta_k^{AWHM} d_k\| \leq \|g_{k+1}\| + |\beta_k^{AWHM}|\|d_k\|$$

$$\|d_{k+1}\| \leq \Gamma + \omega \frac{\|s_k\|}{\lambda_k}$$

$$\|d_{k+1}\| \leq \Gamma + \omega \frac{h}{\lambda^*} = e \Rightarrow \sum_{k \geq 1} \frac{1}{\|d_{k+1}\|^2} = \infty$$

Thus, we find that the global convergence property is fulfilled.

## 3. Result and Discussion

The proposed method was used to solve non-linear convex functions with high dimensions. The comparison was made with the two basic methods that were used in the hybridization process. The comparison was made in terms of execution time and the number of iterations. The efficiency of the method was demonstrated by solving 109 standard problems out of 110 taken from the following references [19, 20]. Figure 1 and figure 2 are prepared according to the Dolan and Moré [21] standard.

The values of the parameters adopted for calculating the step length according to the strong Wolff-Powell conditions: $\boldsymbol{\delta = 10^{-4}, \sigma = 0.9}$

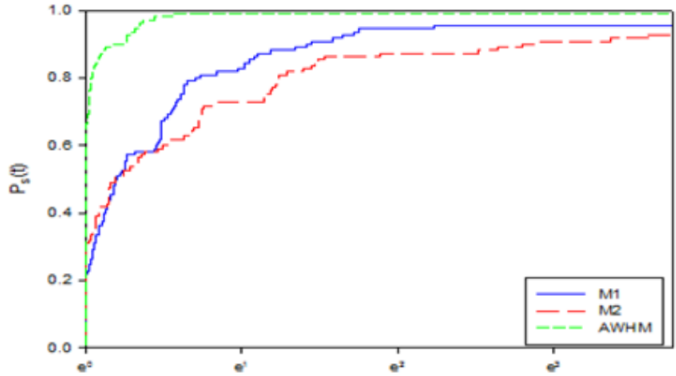

Fig. 1. Comparison of the proposed method with the two basic methods M1 (HRM) and M2 (NHS) in terms of number of iterations

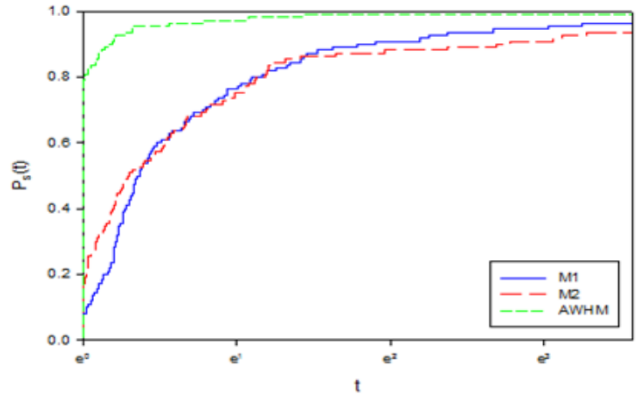

Fig. 2. Comparison of the proposed method with the two basic methods M1 (HRM) and M2 (NHS) in terms of execution time.

*3.1 Application on Arabic medical text*

As an application to this new optimizer, we fine-tuned the BERT model on the Arabic medical dataset [25], in the following, the description of the BERT model and the dataset used.

*3.1.1 Dataset*

The dataset was obtained from three medical volumes (Respiratory System Diseases, Cardiovascular Diseases, and Skin Diseases) issued by the Arabic Encyclopaedia in Syria. The volume of respiratory system diseases volume contains 28 articles with 6691 sentences, where the volume of cardiovascular diseases contains 33 articles with 9464 sentences, and the skin diseases volume contains 22 articles with 5921 sentences [25].
The annotation of entities (disease name, organ name, disease symptoms, and drug name) was done by some Syrian medical students in Syrian universities. Figure 3 shows a screenshot from the dataset.

Fig. 3. Screenshot from the dataset

*3.1.2 BERT Model*
BERT uses the transformer [27] that includes two separate mechanisms (encoder and decoder), where is an encoder reads the text input, the decoder produces a prediction for the task. But since the BERT's goal is to generate a language model, only the encoder is necessary.

The Transformer encoder reads the entire sequence of words at once, not (left-to-right or right-to-left), so it is considered bidirectional.

BERT uses two training strategies, first, one is called "Masked LM" where the model replaces 15% of the words in each sequence with a [MASK] token and attempts to predict the original value of the masked words. The second strategy is called "Next Sentence Prediction (NSP)" where the model receives pairs of sentences and attempts to predict if the second sentence in the pair is subsequent in the original document.

**BERT Fine-tuning**

In the NER task, the model receives a text sequence and is required to mark the various types of entities, as an example in our application (disease name, organ name, disease symptoms, drug name). Using BERT, the model can be trained by feeding the output vector of each token into a classification layer that predicts the NER label. Usually, in the fine-tuning training, most hyper-parameters stay the same as in BERT training. But we have modified the optimizer used in the classification layer that we have added from Adam to our own optimizer. Figure 4 shows the BERT architecture with a fully connected classification layer with the BIO system for the NER task.

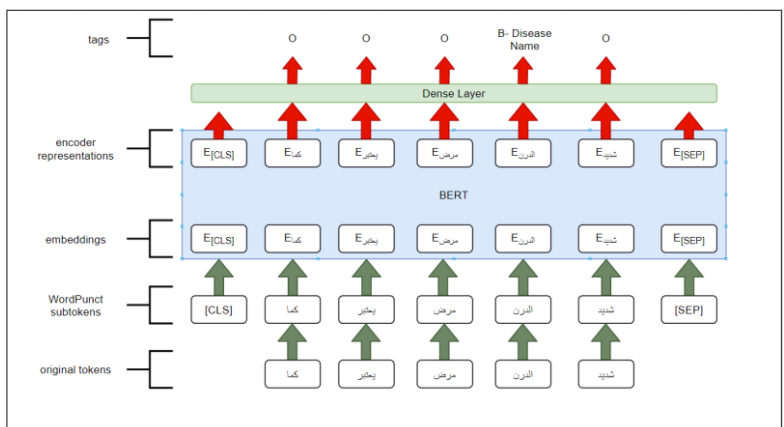

Fig. 4. BERT architecture with a fully connected classification layer with the BIO system for our NER task.

*3.2 Result of the application*

We measure the accuracy of the model by using f1 score:

$$F1 = 2 * \frac{Precision * Recall}{Precision + Recall}$$

Where the precision computed by the following equation:

$$Precision = \frac{TP}{TP + FP}$$

And the recall computed by the following equation:

$$Recall = \frac{TP}{TP + FN}$$

Where TP denotes to true positive, FP to false positive, and FN to false negative. The results were compared with previous results in the following table:

Table 1. F1 score for every entity

| Model Name | Disease Name | Organ Name | Disease Symptoms | Drug Name |
|---|---|---|---|---|
| BERT with Adam optimizer | 87.7102 | 86.3021 | 69.8519 | 77.4986 |
| BERT with new optimizer | 87.9014 | 85.9351 | 68.4125 | 78.0147 |

The results in terms of accuracy show the stability of the proposed method, while the positive results came in terms of execution time.

The same machine was used in the two applications, with 32 batch size, and 4 epochs. The results show an improvement of 18.9 percent over the previous model.

Table 2. Performance Summary

| An example of a column heading | Epochs | Train Time |
|---|---|---|
| BERT with Adam optimizer | 4 | 9.25 hours |
| BERT with new optimizer | 4 | 7.5 hours |

**Conclusion**

Among the most important results that were reached:
- Introducing a new nonlinear gradient method for solving high-dimensional convex functions.
- Study the convergence analysis of the method.
- Demonstrate the effectiveness of the presented method through 110 standard problems of various dimensions and compare it with previous methods.
- Demonstrate the efficiency of the method through direct application to the problem NER in the Arabic medical language, as the results showed the stability of the method and the speed efficiency.

**Reference**

[1] HAN, X., ZHANG, J., CHEN,J., (2017) "**A New Hybrid Conjugate Gradient Algorithm for Unconstrained Optimization**" *Bulletin of the Iranian Mathematical Society* Vol. 43, No. 6, pp. 2067-2084.
[2] HAMODA, M ., MAMAT, M., RIVAIE, M., SALLEH, Z., (2016) "**A Conjugate Gradient Method with Strong Wolfe-Powell Line Search for Unconstrained Optimization**" *Applied Mathematical Sciences*, Vol. 10, No. 15, pp.721 – 734.
[3] SALIH, Y., HAMODA, M., RIVAIE, M., (2018) "**New Hybrid Conjugate Gradient Method with Global Convergence Properties for Unconstrained Optimization**" *Malaysian Journal of Computing and Applied Mathematics*, Vol .1,No1- pp.29-38.
[4] P. Wolfe, (1971) "**Convergence conditions for ascent methods**. II: **some corrections**", *SIAM Review*, vol. 13, no. 2, pp. 185-188.
[5] HESTENES M.R.; STIEFEL E.,(1952) "**Method of Conjugate Gradient for Solving Linear Equations**" *Journal of Research of the National Bureau of Standards*, Vol. 49, No. 6, pp. 409–436.
[6] FLETCHER R. ; REEVES C., (1964) "**Function Minimization by Conjugate Gradients**" *the Computer Journal,* Vol. **7,** No. 2, pp.149-154.
[7] POLAK E.; RIBIÈRE G., (1969) "**Note Sur la Convergence de Méthodes de Directions Conjuguées**" *Revue Françoise d'Informatica et de Recherché Opérationnelle*, Vol.16,No.3, pp. 35–43.
[8] B. T. Polyak, (1969) "**The conjugate gradient method in extreme problems**" *USSR Computational Mathematics and Mathematical Physics*, **9** , 94-112.
[9] GHAZALI, K., SULAIMAN, J., DASRIL, Y., GABDA, D., (2019) "**Newton-MSOR Method for Solving Large-Scale Unconstrained Optimization Problems with an Arrowhead Hessian Matrices**" *Transactions on Science and Technology*, Vol. 6, No. 2-2,pp. 228 – 234.
[10] BERAHAS, A.S., TAK_A_C, M**.,** (2020**) "A robust multi-batch l-bfgs method for machine learning."** *Optimization Methods and Software*, Vol. 35, No.1, pp. 191-219.
[11] BERAHAS, A.S., JAHANI, M., TAK_AC, M., (2019) "**Quasi-Newton Methods for deep learning: Forget the past, just sample**." *arXiv preprint arXiv:* 1901.09997
[12] Khanaiah, Z. and Hmod, G., (2017) "**Novel hybrid algorithm in solving unconstrained optimizations problems**." *Int. J. Novel Res. Phys. Chem. Math.*, *4*(3), pp.36-42
[13] ERWAY, J.B., GRI_N, J., MARCIA, R.F., OMHENI, R., (2019) "**Trust-Region Algorithms for Training Responses: Machine Learning Methods Using Indefinite Hessian Approximations.**" *Optimization Methods and Software* pp. 1-28.
[14] KIMIAEI M., GHADERI S., (2017) "**A New Restarting Adaptive Trust–Region Method for Unconstrained Optimization**." *Journal of the Operations Research Society of China*, 5(4),487–507.
[15] DAWAHDEH, M., et al., (2020) "**A New Spectral Conjugate Gradient Method with Strong Wolfe-Powell Line Search.**" *International Journal of Emerging Trends in Engineering Research*,Vol.8,No.2, pp.391 – 397**.**
[16] ABDULLAHI,I., AHMAD, R., (2017) "**Global Convergence Analysis of A New Hybrid Conjugate Gradient Method for Unconstrained Optimization Problems.**" *Malaysian Journal of Fundamental and Applied Sciences* Vol. 13, No. 2,pp. 40-48**.**
[17] Dai, Y., H., Liao L Z**.,** (2001) "**New conjugacy conditions and related nonlinear conjugate gradient methods.**" *Appl Math Optim*, , Vol.43, pp.87–101**.**
[18] RIVAIE, M., MAMAT, M., ABASHAR, A., (2015) "**A new class of nonlinear conjugate gradient coefficients with exact and inexact line searches.**" *Appl. Math. Comp.*268, pp. 1152-1163.
[19] N. Andrei, (2008**) "An unconstrained optimization test functions collection**." *Advanced Modelling and Optimization*, 10, 147-161.

[20] S. Mishra, (2007) **"Some new test functions for global optimization and performance of repulsive particle swarm method."** *Munich Personal RePEc Archive,* Apr 13.

[21] Dolan, E.D., More, J.J., (2002) "**Benchmarking optimization software with performance profiles**." *Mathematical Programming*, 91 (2), 201-213.

[22] Kingma, D.P. and Ba, J., (2014) "**Adam: A method for stochastic optimization**." *arXiv preprint arXiv:1412.6980*.

[23] Sutskever, I., Martens, J., Dahl, G. and Hinton, G., (2013) "**On the importance of initialization and momentum in deep learning**." In *International conference on machine learning* (pp. 1139-1147). PMLR.

[24] Tieleman, T. and Hinton, G., (2012) "**Lecture 6.5-rmsprop: Divide the gradient by a running average of its recent magnitude."** *COURSERA: Neural networks for machine learning*, *4*(2), pp.26-31.

[25] Hammoud, J., et al., (2020) "**Named entity recognition and information extraction for Arabic medical text",** Multi Conference on Computer Science and Information Systems, MCCSIS 2020 - Proceedings of the International Conference on e-Health 2020, 121-127.

[26] Devlin, J., Chang, M.W., Lee, K. and Toutanova, K., (2018) "**Bert: Pre-training of deep bidirectional transformers for language understanding."** arXiv preprint arXiv:1810.04805.

[27] Vaswani, A., Shazeer, N., Parmar, N., Uszkoreit, J., Jones, L., Gomez, A. N., ... & Polosukhin, I. (2017). "Attention is all you need." *arXiv preprint arXiv:1706.03762*.

[28] Al-Najjar, H., Dabash, M., Al-Issa, A., (2020) **"A Nonlinear Gradient Method for Solving Convex Functions with Large Scale"**, *Al-Baath University Journal for Research and Scientific Studies, Volume (42)*.